\documentclass[journal]{IEEEtran}

\ifCLASSINFOpdf
\else
   \usepackage[dvips]{graphicx}
\fi
\usepackage{url}

\usepackage{graphicx}
\usepackage{amsmath}
\usepackage{amsfonts}

\begin{document}

\title{Novel low-pass filter with adjustable parameters of~exponential-type forgetting}

\author{Ivo Petr\'a\v{s}, \IEEEmembership{Senior Member, IEEE}
\thanks{Manuscript received: October 3, 2022; \\
This work was supported in part by the Slovak Grant Agency for Science under grant VEGA 1/0365/19,  by the Slovak Research and Development Agency under contracts No. APVV-14-0892 and No. APVV-18-0526, and by Army Research Office  under grant No. W911NF-22-1-0264.}
\thanks{Ivo Petr\'a\v{s} is with the Faculty of BERG, Technical University of Ko\v{s}ice, N\v{e}mcovej 3, 042 00, Ko\v{s}ice, Slovak Republic (e-mail: ivo.petras@tuke.sk). \\
}
}

\maketitle

\begin{abstract}
In this paper, a novel form of Gaussian filter, the Mittag-Leffler filter, is presented. This new filter uses a Mittag-Leffler function in the probability density function. Such Mittag-Leffler distribution is used in the convolution kernel of the filter. The filter has three parameters that may adjust the curve shape due to the filter forgetting factor. Illustrative examples present the main advantages of the proposed filter as compared to classical Gaussian filtering techniques. Some implementation notes, together with the Matlab function, are also presented.
\end{abstract}


\IEEEpeerreviewmaketitle

\section{Introduction}

\IEEEPARstart{F}{iltering}
is processing a signal whereby some unwanted components or properties are removed from the signal or some aspects of the signal are suppressed. It often means removing some frequencies or frequency bands from the signal. However, we do not have to use filters exclusively in the frequency domain, and certain frequency components can be removed without having to act in the frequency domain. Filters are widely used, for example, in signal processing in electronics and telecommunications, in radars, control systems sensors, as well as image processing, and computer graphics.
Various forms of filters are used, for instance, the Laplacian filter \cite{Sumiya}, Bayesian filter \cite{Li}, Gaussian filter \cite{Deng}, and so on.

This paper describes a new filter based on the famous Gaussian filter. It is well known that this filter is often used in many areas of signal and image processing for smoothing and noise reduction, e.g., \cite{Hodson, Talmon, Wells, Seddik}. It is a convolutional filter that uses a Gaussian function as a convolution kernel and mathematically adjusts the input signal by convolution with a~Gaussian function. In other words, a Gaussian filter is a filter whose impulse response is a Gaussian function. Among other things, these filters have the important property that they do not overshoot at the input of the step function and, at the same time, minimize the rise and fall time. This behavior is closely related to the Gaussian filter having the minimum possible group delay.

The main contributions of this paper are as follows:
\begin{itemize}
\item generalization of the Gaussian filter to the novel form, based on the Mittag-Leffler distribution function,
\item suggestion of the implementation algorithm for the new Mittag-Leffler filter with adjustable forgetting parameters.
\end{itemize}
The structure of this paper is as follows. Section I briefly
describes the introduction to the problem. Section II presents the essential mathematical tools.
The main results are shown in Section III.
The illustrative examples are presented in Section IV to demonstrate
the benefits of the proposed new filter. Finally, some
concluding remarks are given in Section V.

\section{Preliminaries}

\subsection{Gaussian function and Gaussian distribution}

The Gaussian function, named after Johann Carl Friedrich Gauss is a function that can be expressed in elemental form
\begin{equation}\label{GF}
f(x) = a e^{-\frac{(x-b)^2}{2 c^2}},
\end{equation}
for arbitrary real parameters $a$,\,$ b$, and $c>0$.

Gaussian functions (\ref{GF}) are often used in statistics to represent the probability density function (PDF) of a normal shifted distribution (a.k.a. Gauss distribution) for a real-valued random variable with expected value (or mean) $b=\mu$ and variance $c^2 = \sigma^2$. The general form of its PDF $\phi(x)$ is
\begin{equation}\label{GD}
\phi(x; \sigma) = \frac{1}{\sigma \sqrt{2\pi}} e^{-\frac{1}{2}\left (\frac{x-\mu}{\sigma}\right )^2},
\end{equation}
where the variable $\mu \in \mathbb{R}$ is the mean (or expectation) of the distribution, while the positive variable $\sigma \in \mathbb{R}$ is its standard deviation. The variance of this Gauss distribution is then $\sigma ^{2}$.

The simplest case of the normal  distribution is known as the normal unit distribution or the standard normal distribution. This is a particular case when $\sigma =1$ and $\mu =0$. It means that $x$ has variance, standard deviation of 1, and mean of 0.

Except for the mentioned utilization of the Gaussian function as PDF for normal distribution, we may use it in signal processing to define Gaussian filters and image processing, where a two-dimensional Gaussian filter is used for blurs.

Moreover, the exponential law is the classical approach to studying the dynamics of systems, but there are many systems where dynamics obey a faster or slower law than the exponential law. In that case, the Mittag-Leffler function can best describe such anomalous dynamics changes \cite{Petras}.

\subsection{Mittag-Leffler function and Mittag-Leffler distribution}

The Mittag-Leffler function $E_{\alpha ,\beta}(z)$, named after Magnus Gustaf Mittag-Leffler, is a special function that depends on two parameters, $\alpha$, and $\beta$. It may be expressed by the following series \cite{Podlubny}:
\begin{equation}\label{MLF}
E_{\alpha ,\beta}(z)=\sum_{n=0}^{\infty}\frac{z^n}{\Gamma(\alpha n + \beta)}, \quad \alpha, \beta>0, \, \, z \in \mathbb{C},
\end{equation}
where $\Gamma ( . )$ is the gamma function, for $\beta=1$, we obtain a~one-parameter Mittag-Leffler function $E_{\alpha ,1}(z) \equiv E_{\alpha}(z)$.

The Mittag-Leffler function is sometimes called the queen of the functions \cite{Mainardi}. There are relations between this function and other functions, for instance:
\begin{equation}
E_{1, 1}(z) = e^z, \,\, E_{0^{+},1}(-z^2)=\frac{1}{1+z^2}, \,\, E_{2, 1}(-z^2)=\mbox{cos}(z).
\end{equation}

The Mittag-Leffler function appears naturally in the solution of fractional differential equations, which is essential in fractional calculus theory \cite{Podlubny}. The ordinary and generalized Mittag-Leffler functions interpolate between a purely exponential-law and power-law behavior. It is an important property that may be used in a filter with variable exponential forgetting. However, there are methods where the fractional calculus (fractional-order derivatives/integrals)
can be directly used in filter design and signal processing \cite{Zhang, Chen, Sheng}.

\begin{figure}[!htb]
\centerline{\includegraphics[width=\columnwidth]{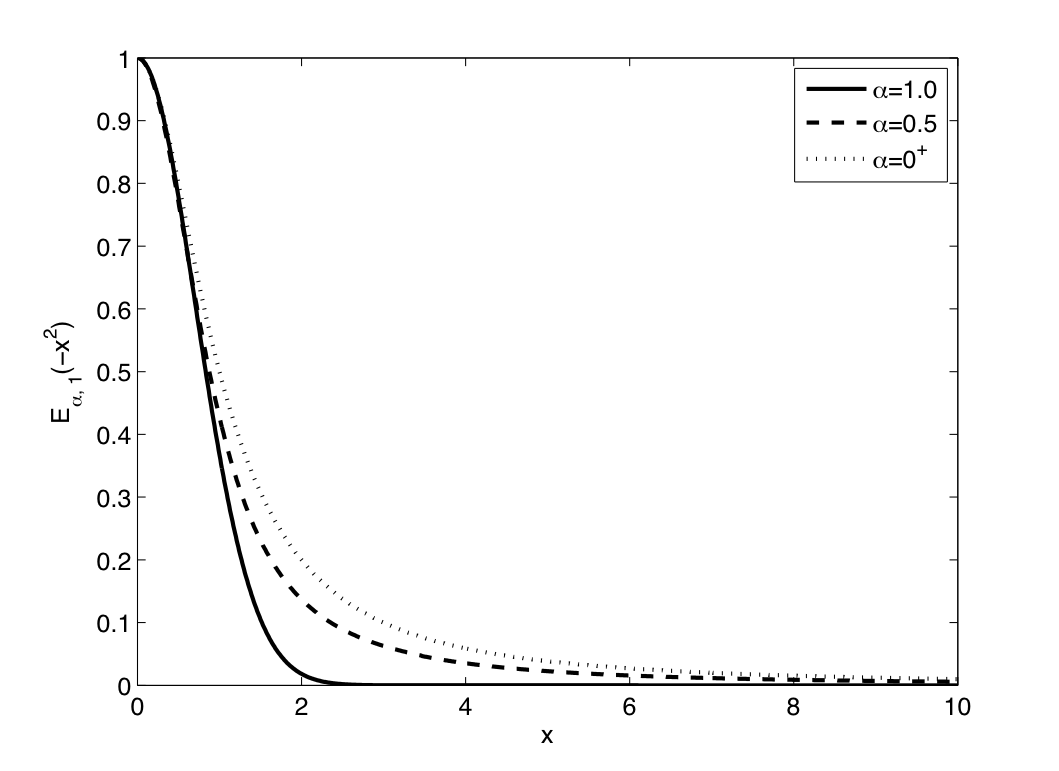}}
\caption{Mittag-Leffler function $E_{\alpha,1}(-x^2)$ for various parameter $\alpha$ within interval $\alpha \in (0; 1]$ and fixed $\beta=1$.}
\label{fig_MLF1}
\end{figure}

\begin{figure}[!htb]
\centerline{\includegraphics[width=\columnwidth]{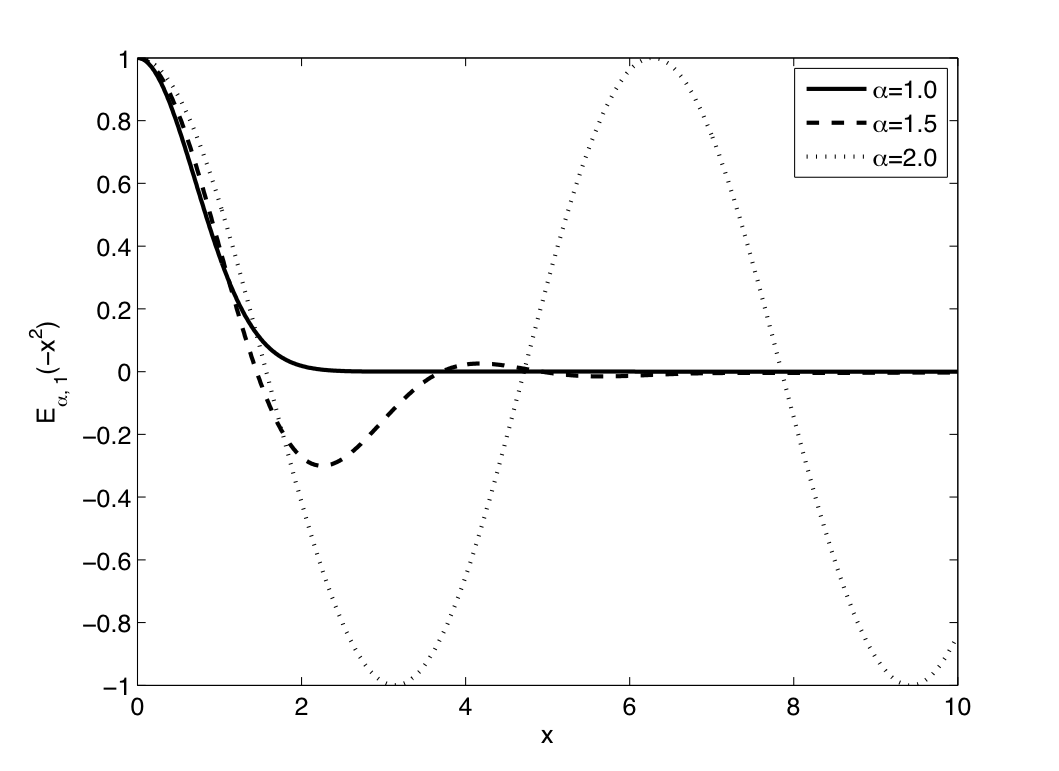}}
\caption{Mittag-Leffler function $E_{\alpha,1}(-x^2)$ for various parameter $\alpha$ within interval $\alpha \in [1; 2]$ and fixed $\beta=1$.}
\label{fig_MLF2}
\end{figure}
Fig.~\ref{fig_MLF1} and Fig.~\ref{fig_MLF2} plot the Mittag-Leffler function $E_{\alpha,1}(-x^2)$ behavior for various values of $\alpha$ and $\beta=1$, respectively.

Furthermore, the Mittag-Leffler distribution with parameter $\alpha$ proposed in \cite{Agahi} can be written by the PDF as:
\begin{equation}\label{MLD}
\phi(x; \sigma, \alpha) = \frac{\sqrt{\sigma}}{\pi} \Gamma \left (\frac{\alpha}{2}\right ) E_{\alpha, \alpha}
\left (-\sigma (x-\mu)^2 \right ),
\end{equation}
where $0 < \alpha \leq 1$,\, $\sigma>0$, and $\mu \in \mathbb{R}$.
Earlier, in 1990, it was proved by Pillai \cite{Pillai} that the Mittag-Leffler distribution with parameter $\alpha$ is
attracted to the stable distribution only with exponent $\alpha$, $0 < \alpha \leq 1$. It is related to the behavior of the Mittag-Leffler function depicted in Fig.~\ref{fig_MLF1}. Hence, function $1 - E_{\alpha, 1}(-x^{\alpha})$ is the cumulative distribution function of a~probability measure on the non-negative real numbers.

Mittag-Leffler distribution with parameters $\alpha$ and $\beta$ can be written by the PDF as \cite{Agahi}:
\begin{equation}
\phi(x; \sigma, \alpha, \beta) =  \frac{\sqrt{\sigma}}{\pi} \Gamma \left (\beta - \frac{\alpha}{2}\right ) E_{\alpha, \beta}
\left (-\sigma (x-\mu)^2 \right ), \label{MLPDF}
\end{equation}
where $\beta \geq \alpha$ in (\ref{MLPDF}). For $\alpha=\beta$, we get the distribution described by (\ref{MLD}). Some other important properties and relations to stochastic processes can be found in \cite{Pillai, Huillet, Agahi, Albrecher}.

However, as we may observe in Fig.~\ref{fig_MLF2}, the behavior of the Mittag-Leffler function for parameter $\alpha>1$ leads to oscillations, and therefore such distribution is related to a~negative probability. Many authors have already observed such behavior, for instance, in 1942 Paul Dirac \cite{Dirac}, and in 1987 Richard Feynman \cite{Feynman}.
Recently, Leonenko and Podlubny have shown that extension of the Monte Carlo approach to fractional differentiation of orders higher than one led to working with signed probabilities that are not necessarily positive \cite{Leonenko}. We adopted this idea for further consideration, and we will use the following form of the Mittag-Leffler distribution with its PDF:
\begin{equation}
\phi(x; \sigma, \alpha, \beta) = \frac{1}{\sigma  \sqrt{2\pi}}  E_{\alpha, \beta}\left (-\frac{ (x-\mu)^2}{2\sigma^2} \right ), \label{MLPDF2}
\end{equation}
where $\sigma$, $\alpha$, and $\beta$ are positive filter parameters, and $\mu$ is the mean value of an independent variable $x$ usually denoted by symbol $\bar x$.

It should also be noted that the Mittag-Leffler distribution is connected to related and other families of distributions \cite{Chakraborty}.

\section{Main results}

\subsection{Formulation}

In general, the objective of the filter is to extract a true signal from
the noisy measured signal
\begin{equation}
y(t)=y_d(t) + y_s(t),
\end{equation}
where $y(t)$ is the observed (measured) signal at the time $t$, $y_d(t)$ is the true, deterministic part of the signal, and $y_s(t)$ is a stationary noise, stochastic (random) part in the signal, which we assume that it has zero mean.

Let us recall the Gaussian low-pass filter in the time domain, which is defined as the convolution of measured (observed) signal $y(t)$, and the Gauss function $\phi(t;\sigma)$ is as follows
\begin{equation}\label{GFi}
y_{Gf}(t)=y(t) \ast \phi(t; \sigma) = \int_{-\infty}^{\infty} y(t-\tau)\phi(\tau; \sigma) d\tau,
\end{equation}
where $y_{Gf}(t)$ is an output from the filter, i.e., filtered signal.

\subsection{Proposed method}

Following the idea of a generalization of the exponential function to the Mittag-Leffler function of two parameters, a~novel filter, let us name it the Mittag-Leffler filter, can be defined as follows
\begin{equation}\label{MLFDef}
y_{MLf}(t)=y(t) \ast \phi(t; \sigma, \alpha, \beta) = \int_{-\infty}^{\infty} y(t-\tau)\phi(\tau; \sigma, \alpha, \beta) d\tau.
\end{equation}

Obviously, we obtain the classical Gaussian filter (\ref{GFi}) for parameters $\alpha=1$ and $\beta=1$.
Moreover, we obtain many filters with adjustable forgetting factors with $\alpha$, $\beta$, and $\sigma$ as tuning knobs. In other words, we can shape the curve of the probability density function. On the other hand, the different distribution shapes mean that we use a different distribution functions, for example, Cauchy distribution, normal distribution, and many others from this family.

Taking into account the Mittag-Leffler distribution (\ref{MLPDF2}) and filter definition (\ref{MLFDef}), we can write the novel filter as:
\begin{equation}
y_{MLf}(t)= \frac{1}{\sigma \sqrt{2\pi}}\int_{-\infty}^{\infty} y(t-\tau) E_{\alpha, \beta}\left (-\frac{ (\tau-\bar \tau)^2}{2\sigma^2} \right ) d\tau. \label{GMLF_eq}
\end{equation}

Such a three-parameter filter is more flexible than the classical one and has more degrees of freedom due to the additional tuning parameters $\alpha$ and $\beta$ with which we can shape the distribution curve and, therefore, exponential forgetting.

\subsection{Implementation notes}

For practical implementation of Gaussian filter in discrete-time domain, we may use methods described, for instance, in
\cite{Rau, Young}.
The Gaussian filter (\ref{GFi}) is not causal, meaning the filter window is symmetric in the time domain. It makes the Gaussian filter physically infeasible because the Gaussian function (\ref{GF}) for $x\in (-\infty ,\infty)$ would theoretically require an infinite window length. For practical implementation, it is reasonable to shorten the filter window and use it directly for narrow windows. However, in some cases, this truncation can cause significant errors. In real-time systems, there is 
a~delay because the incoming samples must fill the filter window before the filter can be applied to the signal being processed. The Gaussian filter kernel in the convolution is continuous. The most common replacement for the continuous kernel is the discrete equivalent sampled Gaussian kernel, represented by sampling points from the continuous Gaussian kernel. Instead of integration operation in convolution, the summation operation over all samples can be used.

It is also well-known that conventional averaging filters based on a moving average, or based on a weighted moving average with exponential forgetting, are not always suitable for their method of assigning weights to older samples of the filtered signal because a frequent request is that the lowest weight was not assigned to the oldest sample, but the sample with a~high proportion of the stochastic component. On the other hand, such filters ensure a quick response to a change to 
a~deterministic or stochastic component by assigning a higher weight to more current components. A requirement is also that the filter algorithm was mathematically and programming relatively simple for use even in digital controllers with limited computing capacity and so that the use of the filter was not limited due to a large number of stochastic values in the measured waveforms.

We may expect similar problems as mentioned above in the Mittag-Leffler filter implementation. Moreover, the Mittag-Leffler function brings some problems with its implementation in real-time applications due to the infinity upper sum limit in the definition. This problem can be circumvented by using the definition of the Mittag-Leffler function in the integral form \cite{Podlubny}. However, there can be a problem with the numerical integration method. These limitations were partially solved, and Podlubny and Kacenak proposed a practical implementation algorithm for the Mittag-Leffler function as a Matlab function \cite{Podlubny-Kacenak}. For further investigation, we use their function:

{
\footnotesize
\begin{verbatim}
function [e]=mlf(alf,bet,c,fi)
%
% MLF -- Mittag-Leffler function.
% MLF(alpha,beta,Z,P) is the Mittag-Leffler function
% E_{alpha,beta}(Z) evaluated with accuracy 10^(-P)
% for each element of Z, alpha and beta are scalars,
% P is integer, Z can be a vector or a 2-dimensional
% array. The output is of the same size as Z.
\end{verbatim}
}

Aforementioned filter (\ref{GMLF_eq}) can be easily implemented as a~Matlab function using function \texttt{mlf(alf,bet,c,fi)}. The Matlab code of the suggested Mittag-Leffler filter (\ref{GMLF_eq}) is the following \cite{Petras_Matlab}:

{
\footnotesize
\begin{verbatim}
function [y_filt] = ML_filter(t,y,sigma,alpha,beta)
%
% Mittag-Leffler filter 
% Inputs: t = independent variable, e.g., time
% y = noisy data to be filtered at points t
% sigma = standard deviation
% alpha, beta = parameters of Mittag-Leffler function
% Output:y_filt = filtered data given in variable y
%
n = length(y);
a = 1/(sqrt(2*pi)*sigma);
dt = diff(t); dt = dt(1);
%
filter=dt*a*mlf(alpha,beta,-0.5*((t-mean(t)).^2)/...
(sigma.^2));
%
y_filt = conv(y, filter, 'same');
Ones2Filter = ones(size(y));
Ones_Filter = conv(Ones2Filter, filter, 'same');
y_filt = y_filt./Ones_Filter;
\end{verbatim}
}

\section{Illustrative examples}

For an illustration of the Mittag-Leffler filter benefits, we use the filter (\ref{GMLF_eq}).
The sampling interval is $0.01$ in both cases.

In the first example, we compare the Gaussian and Mittag-Leffler filters on the noisy signal given by the function
$y_1(t)=e^{-t}\mbox{sin}(3t+1)$ with random noise from a normal distribution.

\begin{figure}[!htb]
\centerline{\includegraphics[width=\columnwidth]{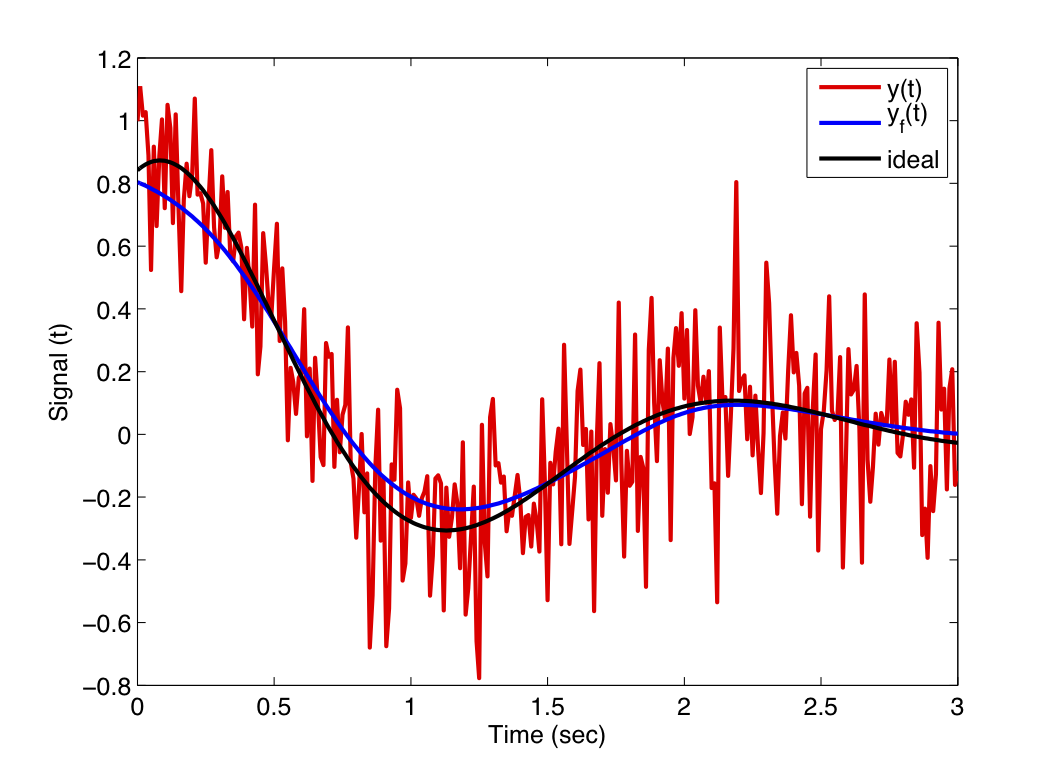}}
\caption{Gaussian filter with parameters $\sigma=0.2$, \,$\alpha=1.00$, \,$\beta=1.00$ applied on noisy test signal $y_1(t)$
and comparison with ideal curve.}
\label{fig_ex1_1}
\end{figure}

\begin{figure}[!htb]
\centerline{\includegraphics[width=\columnwidth]{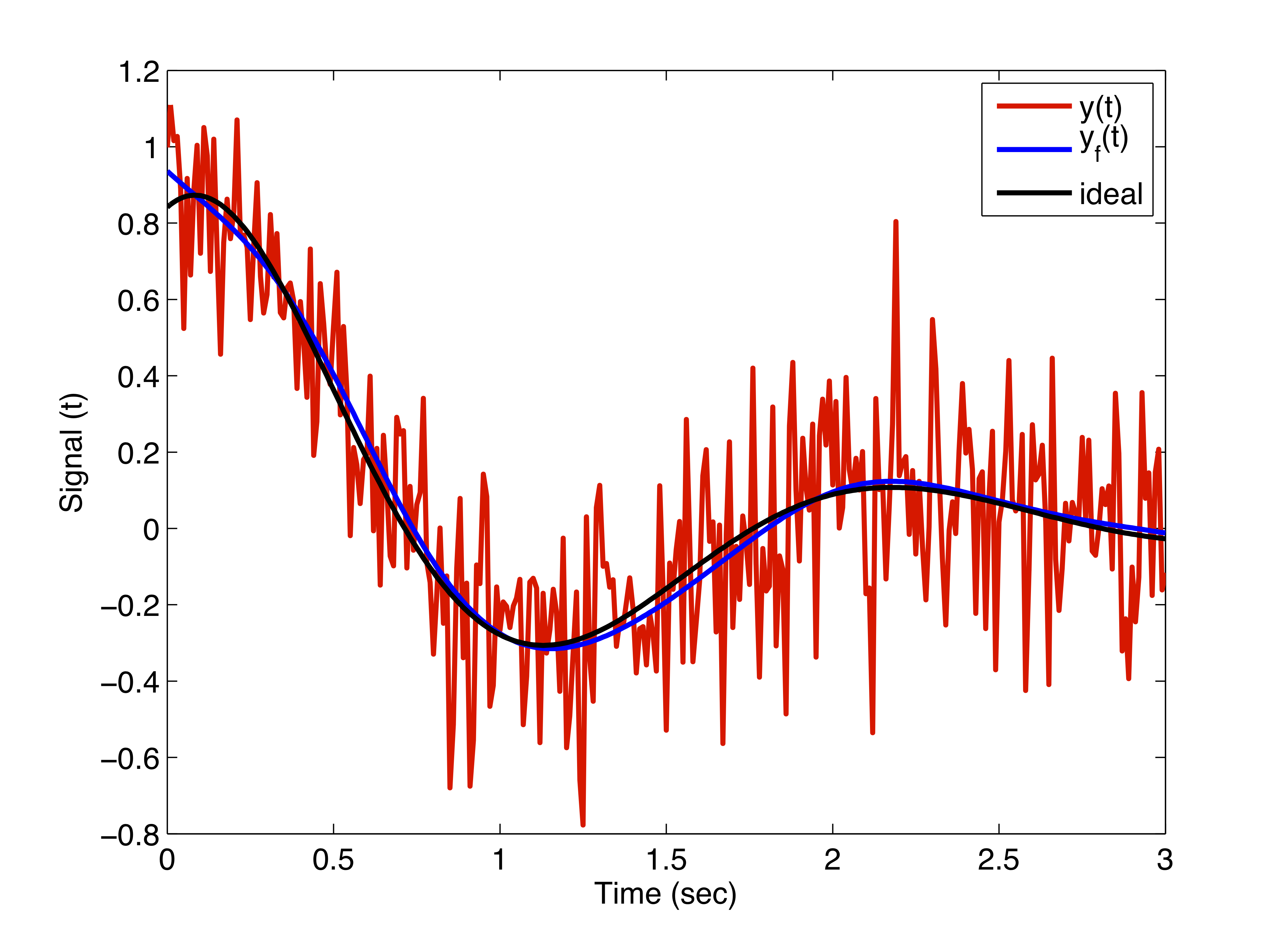}}
\caption{Mittag-Leffler filter with parameters $\sigma=0.2$, \,$\alpha=1.20$, \,$\beta=1.00$ applied on noise test signal 
$y_1(t)$ and comparison with ideal curve.}
\label{fig_ex1_2}
\end{figure}

In Fig.~\ref{fig_ex1_1} is depicted the noisy signal $y_1(t)$, ideal curve without noise, and filtered signal $y_f(t)$ with filter (\ref{GMLF_eq}) and the filter parameters $\sigma=0.2$, \,$\alpha=1.00$, \,$\beta=1.00$. For these parameters $\alpha$ and $\beta$, we obtain the Gaussian filter (\ref{GFi}).

In Fig.~\ref{fig_ex1_2} is shown the noisy signal $y_1(t)$, ideal curve without noise, and filtered signal $y_f(t)$ with filter (\ref{GMLF_eq}) and the filter parameters $\sigma=0.2$, \,$\alpha=1.20$, \,$\beta=1.00$. Thus we get the Mittag-Leffler filter.

In the second example, we compare the Gaussian filter and the Mittag-Leffler filter on the noisy signal given by the function
$y_2(t)=\mbox{sin}(\pi t/0.7) + \mbox{cos}(2\pi t)$ with random noise from a~normal distribution. 

\begin{figure}[!htb]
\centerline{\includegraphics[width=\columnwidth]{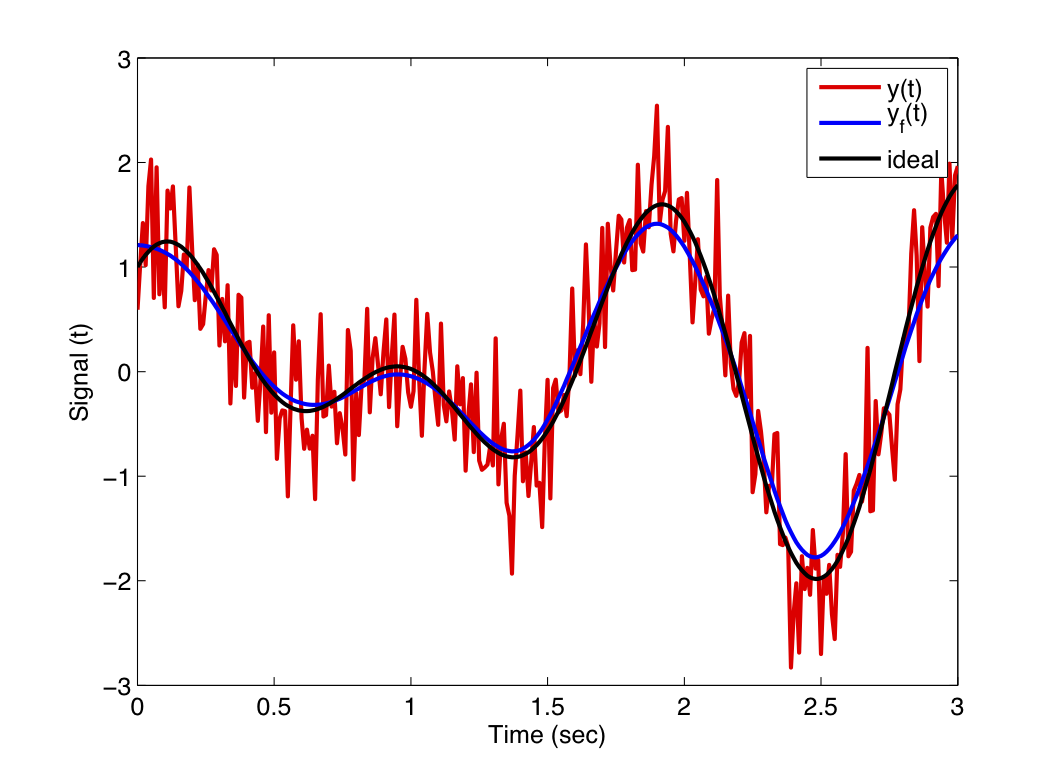}}
\caption{Gaussian filter with parameters $\sigma=0.1$, \,$\alpha=1.00$, \,$\beta=1.00$ applied on noisy test signal $y_2(t)$
and comparison with ideal curve.}
\label{fig_ex2_1}
\end{figure}

\begin{figure}[!htb]
\centerline{\includegraphics[width=\columnwidth]{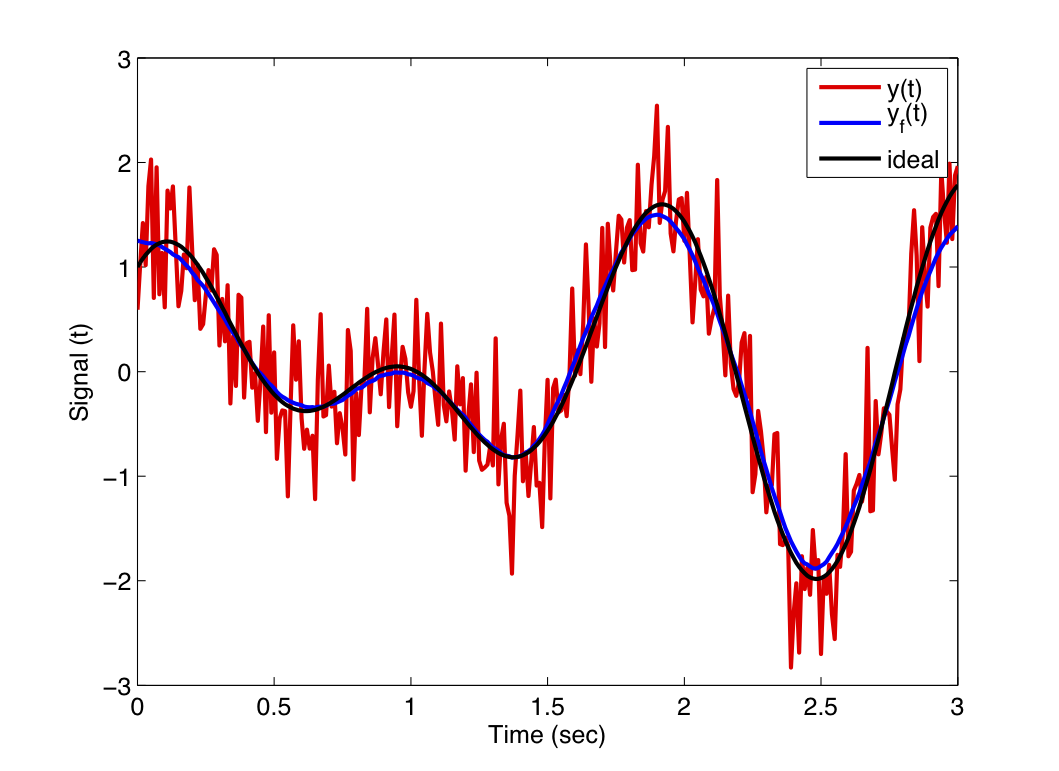}}
\caption{Mittag-Leffler filter with parameters $\sigma=0.1$, \,$\alpha=0.95$, \,$\beta=0.90$ applied on noise test signal 
$y_2(t)$ and comparison with ideal curve.}
\label{fig_ex2_2}
\end{figure}

Fig.~\ref{fig_ex2_1} shows the noisy signal $y_2(t)$, ideal curve without noise, and filtered signal $y_f(t)$ with filter (\ref{GMLF_eq}) and the filter parameters $\sigma=0.1$, \,$\alpha=1.00$, \,$\beta=1.00$. Such parameter values for $\alpha$ and $\beta$ lead to the classical Gaussian filter (\ref{GFi}).

In Fig.~\ref{fig_ex2_2} is presented the noisy signal $y_2(t)$, ideal curve without noise, and filtered signal $y_f(t)$ with filter (\ref{GMLF_eq}) and the filter parameters $\sigma=0.1$, \,$\alpha=0.95$, \,$\beta=0.90$. 
Obviously, for these parameters set, we get the Mittag-Leffler filter.

\begin{table}[!htb]
\caption{Mean Squared Errors Comparison}
\label{table}
\small
\setlength{\tabcolsep}{3pt}
\begin{tabular}{|p{25pt}|p{30pt}|p{125pt}|p{40pt}|}
\hline
Figure & Signal & Filter (\ref{GMLF_eq}) parameters ($\sigma$, $\alpha$, $\beta$)& MSE value \\
\hline
\#\ref{fig_ex1_1} & $y_1(t)$ & $\sigma=0.2$, \,$\alpha=1.00$, \,$\beta=1.00$ & 0.0027 \\
\hline
\#\ref{fig_ex1_2} &$ y_1(t)$ &  $\sigma=0.2$, \,$\alpha=1.20$, \,$\beta=1.00$ & 0.0019 \\
\hline
\#\ref{fig_ex2_1} & $y_2(t)$ & $\sigma=0.1$, \,$\alpha=1.00$, \,$\beta=1.00$ & 0.0233 \\
\hline
\#\ref{fig_ex2_2} & $y_2(t)$ & $\sigma=0.1$, \,$\alpha=0.95$, \,$\beta=0.90$ & 0.0138 \\
\hline
\end{tabular}
\label{tab1}
\end{table}

Table~\ref{tab1} has summarized the comparison of the mean squared errors (MSE) for two test signals and different sets of the filter parameters, which were experimentally found and presented above. As we can see from these first results, the Mittag-Leffler filter gives better results than the Gaussian filter due to more adjustable parameters. It means we can shape the distribution function's curve and thus the exponential-type forgetting factor.

\section{Conclusion}

This paper develops a novel low-pass filter based on the Mittag-Leffler function for signal processing applications.
The proposed filter has three adjustable parameters and is more flexible than a classical Gaussian filter. For a specific set of filter parameters, the Gaussian filter is a particular case of the new Mittag-Leffler filter. A~Matlab function of the suggested Mittag-Leffler filter was created as well. Two illustrative examples present the benefits of the Mittag-Leffler filter by comparing the values of the MSE with the Gaussian filter.

Techniques used in this paper allow us to extend the Mittag-Leffler filter into two dimensions for image processing.


\end{document}